\documentclass[reqno]{conm-p-l}

\usepackage{amscd,amssymb}

\usepackage{fullpage}

\numberwithin{equation}{subsection}
\newtheorem{Theorem}[equation]{Theorem}
\newtheorem{Lemma}[equation]{Lemma}
\newtheorem{Proposition}[equation]{Proposition}
\newtheorem{Corollary}[equation]{Corollary}

\theoremstyle{definition}

\newtheorem{Definition}[equation]{Definition}

\theoremstyle{remark}

\newtheorem{Remark}[equation]{Remark}
\newtheorem{Example}[equation]{Example}

\newcommand{\uln}[1]{\underline{#1}}
\newcommand{\eqdef}{\overset{\text{def}}{=}}
\newcommand{\psb}[1]{[ \! [#1] \! ]}
\newcommand{\lsb}[1]{( \! (#1) \! )}
\newcommand{\tensor}[1]{\underset{#1}{\otimes}}
\newcommand{\fs}[1]{+_{#1}}
\newcommand{\fm}[1]{-_{#1}}

\newcommand{\slot}{\,-\,}
\renewcommand{\:}{\colon}
\newcommand{\xra}[1]{\xrightarrow{#1}}

\newcommand{\CP}{\mathbb{C} \! P^{\infty}}
\newcommand{\Z}{\mathbb{Z}}
\newcommand{\Q}{\mathbb{Q}}
\newcommand{\R}{\mathbb{R}}
\newcommand{\Zp}{\mathbb{Z}_{p}}
\newcommand{\Fp}{\mathbb{F}_{p}}
\newcommand{\Qp}{\mathbb{Q}_{p}}

\newcommand{\C}{\mathbb{C}}

\DeclareMathOperator{\Hom}{Hom}
\DeclareMathOperator{\Ext}{Ext}
\DeclareMathOperator{\spec}{spec}

\DeclareMathOperator*{\colim}{colim}

\DeclareMathOperator{\uHom}{\uln{Hom}}
\DeclareMathOperator{\uExt}{\uln{Ext}}

\newcommand{\hot}{\hat{\otimes}}
\newcommand{\haf}{\hat{\mathbb{A}}^{1}}
\newcommand{\m}[1]{\mathbf{m}_{#1}}

\newcommand{\etale}{\mathrm{et}}
\newcommand{\G}{\mathbb{G}}
\newcommand{\Gm}{\G_m}
\newcommand{\rGm}{\G_m}
\newcommand{\fml}[1]{#1^{0}}         
\newcommand{\ptorsion}{[p^{\infty}]}
\newcommand{\torsion}{_{\mathrm{tors}}}
\newcommand{\Tate}[1]{\mathrm{Tate} (#1)}
\DeclareMathOperator{\LTLifts}{Lifts}

\newcommand{\LT}[1]{\LTLifts_{#1}}
\newcommand{\bF}{\mathbf{F}}
\newcommand{\ZpmZ}{\Z[\tfrac{1}{p}]/\Z}
\newcommand{\epsuniv}{\epsilon^{\mathrm{univ}}}
\newcommand{\padnrm}[1]{\| #1 \|}

\newcommand{\point}{*}
\newcommand{\pr}[1]{p^{#1}}

\newcommand{\TX}{TX}
\newcommand{\LX}{\mathcal{L}X}
\newcommand{\OLX}{\widetilde{\mathcal{L}X}}
\newcommand{\T}{\mathbb{T}}
\newcommand{\smooth}{C^{\infty}}

\newcommand{\norm}{\nu}

\newcommand{\Bo}[2]{{#2}_{#1}}
\newcommand{\triv}[1]{\underline{#1}}

\newcommand{\V}{\mathcal{V}}
\newcommand{\IV}[1]{\mathcal{I}#1}
\newcommand{\In}[2]{I_{#1}#2}

\renewcommand{\P}{\mathbb{P}}
\newcommand{\rc}[1]{\widetilde{#1}}
\newcommand{\eusual}{e_{\text{usual}}}
\newcommand{\alphausual}{\alpha_{\text{usual}}}
\newcommand{\ET}{E_{\T}}
\newcommand{\hET}{\hat E_{\T}}
\newcommand{\eT}{e_{\T}}
\newcommand{\CH}{\hat{H}}
\newcommand{\KT}{K_{\T}}
\newcommand{\borel}[1]{#1_{\mathrm{Borel}}}
\newcommand{\CK}{\hat{K}}
\newcommand{\sinf}{\Sigma^{\infty}}
\newcommand{\gpof}[1]{P_{#1}}
\newcommand{\chars}{\T^{*}}
\newcommand{\odr}[1]{\C (#1)}
\newcommand{\qf}{\hat{q}}
\newcommand{\qm}{q}

\newcommand{\GammaF}{\fml{\Gamma}}

\newcommand{\hT}{\hat{\T}}
\newcommand{\J}[1]{\Bo{\T}{\odr{#1}}}

\begin{document}


\title[Riemann-Roch and the free loop space]{A renormalized
Riemann-Roch formula and the Thom isomorphism for the free loop space}

\author{Matthew Ando}
\address{The University of Illinois at Urbana-Champaign}
\email{mando@math.uiuc.edu}

\author{Jack Morava}
\address{The Johns Hopkins University}
\email{jack@math.jhu.edu}



\begin{abstract}
Let $E$ be a circle-equivariant complex-orientable cohomology theory.
We show that the fixed-point formula applied to the free loop space of
a manifold $X$ can be understood as a Riemann-Roch formula for the
quotient of the formal group of $E$ by a free cyclic subgroup.   The
quotient is not representable, but (locally at $p$) its $p$-torsion
subgroup is, by a $p$-divisible group of height one greater than the
formal group of $E$.  
\end{abstract}

\thanks{Both authors were supported by the NSF}

\date{January 2000}

\maketitle

\begin{quote}
\emph{I believe in the fundamental interconnectedness of all things.}

\hfill ---Dirk Gently \cite{Adams:DirkGently}
\end{quote}

\section{Introduction}

Let $\T$ denote the circle group, and, if $X$ is a compact smooth
manifold, let $\LX\eqdef \smooth (\T,X)$ denote its free loop space.
The group $\T$ acts on $\LX$, and the fixed point manifold is again
$X$, considered as the subspace of constant loops.  In the 1980's,
Witten showed that the fixed-point formula in ordinary equivariant
cohomology, applied to the free loop space $\LX$ of a spin manifold
$X$, yields the index of the Dirac operator (i.e. the $\hat{A}$-genus)
of $X$---a fundamentally $K$-theoretic quantity
\cite{MFA:CircularSymmetry}.  He also applied the fixed-point theorem
in equivariant $K$-theory to a Dirac-like operator on $\LX$ to obtain
the elliptic genus and ``Witten genus'' of $X$
\cite{Witten:Dir}---quantities associated with elliptic cohomology.

Among homotopy theorists, these developments generated
considerable excitement.  The chromatic program organizes the
structure of finite stable homotopy types, locally at a prime $p$,
into layers indexed by nonnegative integers.  The $n$th layer
is detected by a family of cohomology theories $\mathcal{E}_{n}$;
rational cohomology, $K$-theory, and elliptic cohomology are detecting
theories for the first three layers
\cite{Morava:NoetherianLocalizations,DHS,HS}.

The geometry and analysis related to rational cohomology and
$K$-theory are reasonably well-understood, but for $n\geq 2$ and for
elliptic cohomology in particular, very little is known.
Witten's work provides a major
suggestion: for $n=1$ and $n=2$ his analysis gives a correspondence
\begin{equation} \label{eq:correspondence}
\begin{array}{c}
\text{ analysis underlying $\mathcal{E}_{n}$} \\
\text{applied to $X$}
\end{array} \leftrightarrow
\begin{array}{c}
\text{ analysis underlying $\mathcal{E}_{n-1}$} \\
\text{applied to $\LX$.}
\end{array}
\end{equation}

This paper represents our attempt to understand why Witten's procedure
appears to connect the chromatic layers in the manner of
\eqref{eq:correspondence}.  To do this we consider very generally the
fixed-point formula attached to a complex-oriented theory $E$ with
formal group law $F$.  We recall that for $n>0$, such a theory detects
chromatic layer $n$ if the formal group law $F$ has height $n$.

Our first result is that the fixed-point formula of a suitable
equivariant extension of $E$ (Borel cohomology is fine, as is the
usual equivariant $K$-theory) applied to the free loop
space yields a formula which is \emph{identical} to the Riemann-Roch
formula for the quotient $F/ (\qf)$ of the formal group law $F$ by a
free cyclic subgroup $(\qf)$ (compare formulae
\eqref{eq:WittenFormulaRenorm} and \eqref{eq:RR-g}).

The  quotient $F/ (\qf)$ is not a formal group, so to understand its
structure, we work $p$-locally and study its $p$-torsion subgroup $F/
(\qf)\ptorsion$.  We construct a group $\Tate{F}$ with a canonical map
\[
    \Tate{F}\to F/ (\qf),
\]
which induces an isomorphism of torsion subgroups in a suitable
setting.   Our second result is that the group $\Tate{F}\ptorsion
$ is a $p$-divisible group, fitting into an extension
\[
   F\ptorsion \rightarrow \Tate{F} \rightarrow \Qp/\Zp
\]
of $p$-divisible groups.  If the height of $F$ is $n$, then the height
of $\Tate{F}\ptorsion$ is $n+1$, but its \'etale quotient has height
$1$.  In a sense we make precise in \S \ref{sec:universalprops}, it is the
universal such extension.

Thus the fixed-point formula on the free loop space interpolates
between the chromatic layers in the same way that $p$-divisible groups
of height $n+1$ with \'etale quotient of height $1$ interpolate
between formal groups of height $n$ and formal groups of height $n+1$.
This is discussed in more detail, from the homotopy-theoretic point of
view, in our earlier paper \cite{AMS:Tate} with Hal Sadofsky; this paper
is a kind of continuation, concerned with analytic aspects of these
phenomena. We show that Witten's construction in rational cohomology
produces $K$-theoretic genera because of the exponential exact sequence
\begin{equation}\label{eq:gm-exact-seq}
   0 \rightarrow \Z \rightarrow \C \rightarrow \C^{\times} \rightarrow 1
\end{equation}
expressing the multiplicative group ($K$-theory) as the quotient of
the additive group (ordinary cohomology) by a free cyclic subgroup;
while his work in $K$-theory produces elliptic genera because of the
exact sequence
\begin{equation}\label{eq:tate-exact-seq}
   0 \rightarrow q^{\Z} \rightarrow \C^{\times} \rightarrow
\C^{\times}/q^{\Z} \rightarrow 1
\end{equation}
(where $q$ is a complex number with $|q|<1$), expressing the Tate elliptic
curve $\C^{\times}/q^{\Z}$ as the quotient of the multiplicative group
by a free cyclic subgroup.

These analytic quotients have already been put to good use in
equivariant topology.
Grojnowski constructs from equivariant ordinary cohomology a complex
$\T$-equivariant elliptic cohomology
using the elliptic curve $\C/\Lambda$ which is the quotient of the
complex plane by a lattice; and Rosu uses Grojnowski's functor to give a
striking conceptual proof of the rigidity of the elliptic genus.
Grojnowski's ideas applied to the multiplicative sequence
\eqref{eq:tate-exact-seq} give a construction of complex
$\T$-equivariant elliptic cohomology based on equivariant $K$-theory;
details will appear elsewhere. Completing this circle, Rosu has used
the quotient \eqref{eq:gm-exact-seq} to give a construction of complex
equivariant $K$-theory \cite{Grojnowski:ell,Rosu:rigidity,KnutsonRosu}.

Several of the formulae in this paper involve formal infinite
products; see for example \eqref{eq:WittenFormulaRenorm} and
\eqref{eq:RR-g}.  On the fixed-point formula side, the
source of these is the Euler class of the normal bundle $\nu$ of $X$ in
$\LX$ \eqref{eq:euler-normal}.  From this point of view, the problem
is that the bundle $\nu$ is infinite-dimensional, so it does not have
a Thom spectrum in the usual sense. However, $\nu$ has a highly nontrivial
circle action, which defines a locally finite-dimensional filtration
by eigenspaces. Following the program sketched in \cite{CJS:Floer},
 we construct from this filtration a
Thom \emph{pro-spectrum}, whose Thom class is the infinite product.

In the particular cases of the additive and multiplicative formal
groups ($n=1,2$ above), one can also control the infinite products by
replacing them with products which converge to holomorphic functions
on $\C$; this construction of elliptic functions goes back to Eisenstein.
We are grateful to Kapranov for pointing out
to us that Eistenstein considered the the analogous problem for
$n>2$.   In \cite{Eisenstein:BEAT} he described the difficulty of
interpreting such infinite products.  He went on to hint that he
perceived a useful approach, and concluded the following.

\begin{quote}
Die Functionen, zu welchen man auf diesen Wege gef\"uhrt wird,
scheinen sehr merkw\"urdige Eigenschaften zu besitzen; sie er\"offnen
ein Feld, auf dem sich Stoff zu den reichhaltigsten Untersuchungen
darbietet, und welches der eigentliche Grund und Boden zu sein
scheint, auf welchem die schwierigsten Theile der Analysis und
Zahlentheorie ineinander greifen. 
\end{quote}

\subsection{Formal group schemes}

In this paper (especially in section \ref{sec:GeneralizedTateCurve})
we shall consider formal
schemes in the sense of
\cite{Strickland:FSFG,Demazure:pdivisiblegps}.  A \emph{formal scheme}
is a filtered colimit of affine schemes.  For example the
``formal line''
\[
   \haf \eqdef \colim_{n} \spec \Z[x]/x^{n}
\]
is a formal scheme.  Note that an affine scheme is a formal scheme in a
trivial way. An important feature of this category which we shall
use is that it has finite products.  For example,
\[
   \haf\times \haf =
   \colim \spec\bigl( \Z[x]/ (x^{n}) \otimes \Z[y]/ (y^{m})\bigr).
\]
In particular a \emph{formal
group scheme} means an abelian group  in the category of formal
schemes.  A formal group scheme whose underlying formal scheme
is isomorphic to the formal scheme $\haf$ is called a commutative
one-dimensional formal Lie group.  We shall simply call it a \emph{formal
group}.

The first reason for considering formal schemes is that formal groups
are not quite groups in the category of affine schemes, because a
group law
\[
     F (s,t) = s + t + \dots  \in R\psb{s,t}
\]
over a ring $R$ gives a diagonal
\[
    R\psb{s} \xra{} R\psb{s,t} \cong R\psb{s}\hot R\psb{t}
\]
only to the completed tensor product.

The second reason for considering formal schemes is that, if $G$ is an
affine group scheme, then its torsion subgroup $G\torsion$ is a formal
scheme (the colimit of the affine schemes $G[N]$ of torsion of order
$N$), but not in general a scheme.

If $X$ is a formal scheme over $R$, and $S$ is an $R$-algebra, then
$X_{S}$ will denote the resulting formal scheme over $S$.

\section{The umkehr homomorphism and an ungraded analogue}
\label{sec:umkehr}

\subsection{} \label{sec:umkehr-setup}

Let $E$ be a complex-oriented multiplicative cohomology theory with
formal group law $F$, and let   $h\: X \to Y$ be a proper complex-oriented map
of smooth finite-dimensional connected manifolds, of fiber dimension
$d=\dim X - \dim Y$.  The Pontrjagin-Thom collapse associates to these
data an ``umkehr'' homomorphism
\cite{Quillen:Ele}
\[
      h_{!} : E^{*}(X) \rightarrow E^{*-d}(Y).
\]
We will be concerned with similar homomorphisms in certain
infinite-dimensional contexts.  In order to do so, we
systematically eliminate the shift of $-d$ in the degree by restricting
our attention to \emph{even periodic} cohomology theories $E$.  The
examples show \eqref{sec:examples} that this amounts to measuring
quantities relative to the vacuum.

\subsection{Even periodic ring theories}

Let $E$ be a cohomology theory.
If $X$ is a space,  then $E^{*}(X)$ will denote its \emph{unreduced}
cohomology; if $A$ is a spectrum, then $E^{*} (A)$ will denote its
cohomology in the usual sense.   These notations are related by the
isomorphism $E^{*} (X)\cong E^{*} (\sinf X_{+})$, where $X_{+}$
denotes the union of $X$ and a disjoint basepoint.  The reduced
cohomology of $X$ will be denoted $\rc{E} (X)$.  Let $\point$ denote
the one-point space.

A cohomology theory $E$ with commutative multiplication is \emph{even}
if $E^{\text{odd}}(\point)=0$.  It is \emph{periodic} if
$E^{2}(\point)$ contains a unit of $E^{*}(\point)$.
If $E$ is an even periodic theory, then we write $E(X)$ for $E^{0}(X)$
and $E$ for $E^{0}(\point)$.  We sometimes write $X_{E}= \spec
E(X)$ for the spectrum, in the sense of commutative algebra, of the
commutative ring $E(X)$.

A space $X$ is \emph{even} if $H_{*}(X)$ is a free abelian group,
concentrated in even degrees.  In that case the natural map
\begin{equation} \label{eq:filtration}
    \colim F_{E} \to X_{E},
\end{equation}
where $F$ is the filtered system of maps of finite CW complexes to
$X$, is an isomorphism. This gives $X_{E}$ the structure
of a formal scheme.  The functor $X\mapsto X_{E}$
from even spaces to formal schemes over $E$ preserves finite products
and coproducts: if $X$ and $Y$ are  two even spaces, then
\[
(X\times Y)_{E} \cong X_{E} \times Y_{E}  \cong \spec E(X)\hot E(Y).
\]
Here $\hot$ refers to the completion of the tensor
product with respect to the topology defined by the filtrations of $E(X)$
and $E(Y)$.

\subsection{Orientations and coordinates}

Let $P\eqdef \CP$ be the classifying space for complex line bundles.
Let $m\:  P\times P \to  P$ be the map classifying the tensor product
of line bundles.  It induces a map
\[
\gpof{E}\times \gpof{E} \xra{m_{E}} \gpof{E},
\]
which makes $\gpof{E}$ a formal group scheme over $E$.
Of course it is a formal group: let
$i\: S^{2}\to  P$ be the map classifying the Hopf bundle.
A choice of element $x\in \rc{E}( P)$ such that $v=i^{*}x \in \rc{E}(S^{2})
\cong E^{-2} (\point)$ is a unit is called a \emph{coordinate} on
$\gpof{E}$.
%
There is then an isomorphism
\[
      E( P) \cong E\psb{x},
\]
which determines a formal group \emph{law} $F$ over $E$ by the formula
\[
    F (x,y) = m^{*}x \in E ( P\times P) \cong E\psb{x,y}.
\]
Any even-periodic cohomology theory $E$ is complex-orientable.
An \emph{orientation} on $E$ is a multiplicative natural
transformation
\[
    MU \to E.
\]
These correspond bijectively
with elements $u\in \rc{E}^{2} (P)$ such that
\begin{equation} \label{eq:normalize-orientation}
i^{*}u=\Sigma^2(1),
\end{equation}
where $\Sigma$ is the suspension isomorphism \cite{Adams:BlueBook}.
A coordinate $x$ thus determines an orientation $u=v^{-1}x$.
\begin{Definition} \label{def:oriented-theory}
We shall use the notation $(E,x,F)$ to denote an even periodic
cohomology theory $E$ with coordinate $x$ and group law $F$.  We shall
call such a triple a \emph{parametrized} theory.
\end{Definition}

\subsection{Thom isomorphism}

An orientation $u\in \rc{E}^{2} (P)$ gives the usual Thom classes and
characteristic  classes for complex vector bundles.
If $k$ is an integer, let $\triv{k}$ denote
the trivial complex vector bundle of rank $k$.  If $X$ is a connected
space and $V$ is a complex vector bundle of rank $d$ over $X$, then
we write
\[
X^{V} \eqdef \sinf (\P (V \oplus \triv{1})/\P (V))
\]
for the suspension spectrum of its Thom space, with bottom cell in
degree $2d$.
We write $\alphausual^{V}$ for the Thom isomorphism
\[
     \alphausual^{V}\: E^{*} (X) \cong E^{*+2d} (X^{V}).
\]

In the same way, a coordinate $x\in \rc{E} (P)$ gives rise to a Thom
isomorphism
\[
\alpha^{V}\: E (X)\cong E (X^{V}).
\]
If $v=i^{*}x\in \rc{E} (S^{2})$ is the associated orientation, the
isomorphisms $\alphausual$ and $\alpha$ are related by the formula
\[
     \alpha^{V} = v^{\text{rank }V} \alphausual^{V}.
\]

\begin{Remark} \label{rem:thom-isos}
One effect of condition \eqref{eq:normalize-orientation} is that
$\alphausual^{\triv{d}}$ coincides with the suspension isomorphism
\[
    \alphausual^{\triv{d}} = \Sigma^{2d}\:
    E^{*} (X)\cong E^{*+2d} (X^{\triv{d}}).
\]
The Thom isomorphism $\alpha$ defined by a coordinate chooses $v\in
\rc{E} (S^{2}) \cong E (\point^{\triv{1}})$ as $\alpha^{\triv{1}}$.
Thus $\alphausual$ may be viewed as a composition of Thom isomorphisms
\[
 \alphausual^{V}\: E (X^{\triv{d}}) \xra{(\alpha^{\triv{d}})^{-1}}
                   E (X)            \xra{\alpha^{V}}
                   E (X^{V}).
\]
\end{Remark}

If $\zeta\: \sinf X_{+} \to X^{V}$ denotes the zero section, then we
write
\begin{align*}
    \eusual (V) & \eqdef \zeta^{*}\alphausual^{V} (1) \in E^{2d} (X)  \\
    e (V) & \eqdef \zeta^{*} \alpha^{V} (1) \in E (X)
\end{align*}
for the usual and degree-zero Euler classes of $V$; these are related
by the formula
\[
     e (V)  = v^{\text{rank }V} \eusual (V).
\]

If $U(n)$ denotes the unitary group and $T$ is its maximal
torus of diagonal matrices, then the map
\[
   E(BU(n)) \to E(BT)\cong E ((B\T)^{n})
\]
is the inclusion of the ring of invariants under the action of the
Weyl group $W$.   The coordinate gives an isomorphism
\[
   E(BT) \cong E((B\T)^n) \cong E\psb{x_1,...,x_n},
\]
with $W$ acting as the permutation group $\Sigma_{n}$ on the
$x_{i}$'s.  Thus we can define degree-zero Chern classes $c_i$ in
$E(BU(n))$ by the formula
\begin{equation}\label{eq:chern-splitting}
       \sum_{i=0}^{n} c_i z^{n-i} = \prod_{i=0}^{n} (z + x_i).
\end{equation}
If $F$ is the group law resulting from the coordinate $x$, then we
call the $c_i$ the ``$F$--Chern classes''.

Returning to the map
\[
  X \xrightarrow[]{ h } Y,
\]
as in \eqref{sec:umkehr-setup}, we can now define an umkehr map
\[
  E(X) \xrightarrow[]{ h_{F} } E(Y),
\]
using the degree-zero Thom isomorphism $\alpha$.  We write $F$ to
indicate the dependence on the coordinate.

\begin{Definition} \label{def:genus}
If $X$ is any manifold, we denote by $\pr{X}$ the map
\[
  X \xrightarrow[]{ \pr{X} } \point.
\]
If $E$ is an even periodic theory with group law
$F$, then its $F$--\emph{genus} is the element
$\pr{X}_{F}(1)$ of  $E$.
\end{Definition}

\subsection{The Riemann-Roch formula}

The Riemann-Roch formula compares the umkehr homomorphisms $h_{F}$
and $h_{G}$ of two coordinates with formal group laws $F$
and $G$, related by an isomorphism
\[
    \theta \: F \rightarrow G.
\]
The book of Dyer \cite{Dyer:Cohomology} is a standard reference.

\begin{Proposition} \label{Th:RiemannRoch}
If $h \: X \to Y$  is a proper complex-oriented map of fiber dimension
$2d$, then
\begin{equation}\label{eq:RiemannRoch}
  h_{G}(u) = h_{F}\left[u \cdot \prod_{j=1}^{d}
  \frac{x_{j}}{\theta(x_{j})}\right] \;,
\end{equation}
where the $x_{i}$ are the terms in the factorization
\[
 z^d  + c_{1} z^{d-1} +  \ldots + c_{d} = \prod_{j=1}^{d} (z + x_{i})
\]
of the total $F$--Chern class of the formal inverse of the normal bundle
of $h$. \qed
\end{Proposition}

\begin{Remark}
Changing the coordinate by a unit $u \in E$ multiplies
the umkehr homomorphism by $u^d$; by such a
renormalization, we can always assume that $\theta$ is a
\emph{strict} isomorphism.
\end{Remark}

\subsection{The fixed-point formula}

\subsubsection{Notations for circle actions}

Let $\T$ denote the circle group $\R/\Z$.  If $X$ is a
$\T$-space then we write
\[
   \Bo{\T}{X} \eqdef E\T \underset{\T}{\times} X
\]
for the Borel construction and $X^\T$ for the fixed-point set.
Let $\chars = \Hom[\T,\C^{\times}]$ be the character group of $\T$;
we will also write $\hT = \chars - \{ 1 \}$ for the set of nontrivial
irreducible representations. For $k\in \chars$, let $\odr{k}$ be
the associated one-dimensional complex representation. There is then
an associated complex line bundle $\J{k}$ over $B\T$.

It is convenient to choose an an isomorphism $\T^{*} \cong \Z$;
this determines, in particular, an isomorphism $B\T \cong \CP$. For
$k \in \Z$ we have $\odr{k} = \odr{1}^{\otimes k}$ and
$\J{k} = \J{1}^{\otimes k}$. If $\qf \in E(B\T)$ is the Euler class
of $\J{1}$, then the Euler class of $\J{k}$ is $[k](\qf)$.

\subsubsection{Equivariant cohomology}

\begin{Definition}
Let $(E,x,F)$ be a parametrized theory.  A $\T$-equivariant cohomology
theory $\ET$ is an \emph{extension} of $(E,x,F)$ if
\begin{enumerate}
\item [(1)] There is a natural transformation
\[
     E (X/\T) \rightarrow \ET (X),
\]
which is an isomorphism if $\T$ acts freely on $X$.  In particular the
coefficient ring $\ET (\point)$ is an algebra over $E (\point)$, and
so it is $2$-periodic.
\item [(2)] There is a natural forgetful transformation
\[
     \ET (X) \to E (X).
\]
If $X$ is a trivial $\T$-space then the composition
\[
    E (X) \rightarrow \ET (X) \rightarrow E (X)
\]
is the identity.
\item [(3)] $\ET$ has Thom classes and so Euler classes for complex
$\T$-vector bundles, which are multiplicative and natural under
pull-back.  If $V/X$ is such a bundle, then we write $\eT (V)\in \ET
(X)$ for its (degree-zero) Euler class.  These are compatible with the
Thom isomorphism in $E$ in the sense that, if the $\T$-action on $V/X$
is trivial, then
\[
   \eT (V) = e (V).
\]
\item [(4)] If $L_{1}$ and $L_{2}$ are complex $\T$-line bundles, then
\[
    \eT (L_{1}\otimes L_{2}) = \eT (L_{1}) \fs{F} \eT (L_{2}).
\]
\end{enumerate}
\end{Definition}

\begin{Definition} \label{def:suitable}
If $\ET$ is equivariantly complex oriented as above, a homomorphism
$\ET \rightarrow \hET$ of multiplicative $\T$-equivariant cohomology
theories is a \emph{suitable} localization if
\begin{enumerate}
\item [(1)]
$\hET(\point)$ is flat over $\ET(\point)$,
\item
When $k \neq 0$, $\eT(\odr{k})$ maps to a unit of $\hET(\point)$, and
\item [(2)]
The fixed-point formula \eqref{eq:fixedpointformula} holds for $\hET$.
\end{enumerate}
\end{Definition}

In order to state the fixed-point formula, we need the
following observation of \cite{AtiyahSegal:IndexII}.

\begin{Lemma} \label{t-le:localization}
Let $\ET$ be a suitable theory.  Let $S$ be a compact manifold with
trivial $\T$--action, and let $V$ be a complex $\T$--vector bundle
over $S$.   If the fixed-point bundle $V^\T$ is zero, then
$\eT(V)$ is a unit of $\ET(S)$.
\end{Lemma}

\begin{proof}
Recall \cite{Segal:EquivariantK} that the natural map
\[
 \bigoplus_{k\in \hT} V(k)\otimes \odr{k} \xra{} V
\]
is an isomorphism, where $V(k)\eqdef \Hom[\odr{k},V]$ is the
evident vector bundle over $S$ with trivial $\T$--action.  By
applying the ordinary splitting principle to $V(k)$, we are
reduced to the case that $V=L\otimes \odr{k}$, where $L$ is a
complex line bundle over $S$ with trivial $\T$-action.
If $\ET$ is suitable then the Euler class of $V$ is
\[
\eT(L\otimes \odr{k}) = e(L) \fs{F} \eT(\odr{k}).
\]
Since $S$ is a compact manifold, $e(L)$ is nilpotent in $E(S)$,
so $\eT(V)$ is a unit of $\ET(S)$ because $\eT(\odr{k})$ is a unit
of $\ET(\point)$.
\end{proof}

Now suppose that $M$ is a compact almost-complex manifold with a
compatible $\T$-action.  Let
\[
       j \: S \rightarrow M
\]
denote the inclusion of the fixed-point set; it is a
complex-oriented equivariant map, with $\T$--equivariant normal
bundle $\norm$.  The fixed-point formula which we require in
Definition \ref{def:suitable} is the equation
\begin{equation}
\label{eq:fixedpointformula}
p^{M}_{F}(u) = p^{S}_{F} \left(\frac{j^{*}u}{\eT(\norm)}\right) \;.
\end{equation}
By Lemma \ref{t-le:localization} $\eT(\norm)$ is a unit of $\ET(S)$,
so this localization theorem is a corollary to the projection formula
\[
j_Fj^*(x) = x \cdot \eT(\norm)
\]
for the umkehr of the inclusion of the fixed-point set.

\begin{Example}\label{ex:borel}
The Borel extension
 \[
 \borel{E}(X) \eqdef E(\Bo{\T}{X})
 \]
of an even periodic ring theory has Thom classes $\eT(V) =
e(\Bo{\T}{V})$ for complex $\T$-vector bundles, and the
localization defined by inverting the multiplicative
subset generated by $\eT(\odr{k}), k \neq 0$ will be suitable.
\end{Example}

\begin{Example} \label{ex:kthy}
Let $\KT$ denote the usual equivariant $K$-theory.  Then $K_{\T} =
\Z[q,q^{-1}]$, where $q$ is the representation $\odr{1}$, considered as a
vector bundle over a point.  The Euler class of a line bundle is
$\eT (L) = 1 - L$, so $\qf = 1 - q$.  The group law is multiplicative:
\begin{equation}\label{eq:gm}
   \Gm (x,y) = x + y - x y.
\end{equation}
We have $[k](\qf) = 1 - q^{k}$, and consequently
\[
\CK(X) \eqdef  K_{\T}(X) \otimes_{K_{\T}} \Z((q))
\]
is suitable.
\end{Example}

\begin{Example} \label{ex:CH}
If
\[
H^{*}_{\T}(X)\eqdef H^{*}(\Bo{\T}{X};\Q[v,v^{-1}])
\]
is Borel cohomology with two-periodic rational coefficients, and
$\qf = e(\J{1})$, then $H_{\T}(\point) \cong \Q\psb{\qf}$,
and $e(\J{k}) = k\qf.$  The rational Tate cohomology
\[
\CH^{*}(X)  \eqdef H^{*}_{\T}(X)[\qf^{-1}]
\]
is suitable.
\end{Example}

\begin{Example}
More generally, two--periodic $\hat K(n)_{\T}$ (with $n$ finite
positive) is suitable: if $[p]_{K(n)}(X) = X^{p^n}$ and $k = k_0p^s$
with $(k_0,p) = 1$ then
\[
[k](\qf) = [k_0](\qf^{p^{ns}}) = k_0 \qf^{p^{ns}} + \dots \in \Fp((\qf))
\]
has invertible leading term. Integral lifts of $K(n)$ behave similarly;
the Cohen ring \cite{AMS:Tate} of $\Fp((\qf))$ defines a completion of
the Borel-Tate localization.
\end{Example}

These coefficient rings have natural topologies, which are relevant
to the convergence of infinite products in Corollary
\ref{t-co:procomplex-coh-thom}.

\section{Application to the free loop space}

Let $\hET$ be a suitable localization \eqref{def:suitable} of an
equivariantly complex oriented cohomology theory, let $X$ be a
compact complex-oriented manifold, and let $\LX$ be its free
loop space. Since $\LX$ is not finite-dimensional, the existence
of an umkehr homomorphism $p_{F}^{\LX}$ is not clear.
However, $\T$ acts on $\LX$ by rotations with fixed set $X$ of
constant loops, and Witten discovered that the fixed-point
formula \eqref{eq:fixedpointformula} for the $F$-genus
$p_{F}^{\LX}(1)$ of $\LX$ continues to yield interesting
formulae. In this section we review his calculation.

\subsection{The normal bundle to the constant loops and its Euler
class}

One approximates the space $\smooth (S^{1},\C)$
by the sub-vector space of Laurent polynomials
\[
  \C[\chars] \cong \bigoplus_{k \in \T^{*}} \odr{k} \hookrightarrow
\smooth (S^{1},\C).
\]
The tangent space of $\LX$ is
$\smooth (S^{1},\TX)$.  If $p\in X$ is considered as a
constant loop, then the tangent space to $\LX$ at $p$ is the
$\T$--space $T\LX_{p}\cong \smooth (S^{1},\TX_{p})$.  It is a
$\T$-bundle with a Laurent polynomial approximation
\[
     \TX_p\otimes \C[\chars]
\]
Thus the normal bundle $\norm$ of the inclusion of $X$ in $\LX$ has
approximation
\begin{equation}\label{eq:normalbundle}
     \norm \simeq \bigoplus_{k\in \hT} \TX\otimes \odr{k}.
\end{equation}

If
\[
  z^d  + c_{1}z^{d-1} + \ldots + c_{d} = \prod_{j=1}^{d} (z + x_{i})
\]
is the formal factorization  of the total
$F$-Chern class  of $\TX$, then
\begin{equation} \label{eq:euler-normal}
  \eT(\norm) = \prod_{j = 1}^{d} \prod_{k \neq 0} (x_{j} \fs{F}
  [k]_{F}(\qf)),
\end{equation}
where $\qf = \eT (\odr{1})$.

\subsection{The fixed point formula}
Applying (\ref{eq:fixedpointformula}) to the inclusion
\[
  X \xrightarrow[]{} \LX
\]
yields the formula
\begin{equation} \label{eq:WittenFormula}
  p^{\LX}_{F}(1) = \pr{X}_{F}\left[\prod_{j = 1}^{d} \prod_{k \neq 0}
  \frac{1}{x_{j} \fs{F} [k]_{F}(\qf)} \right].
\end{equation}

Equation \eqref{eq:WittenFormula} requires some interpretive
legerdemain.  For example, the leading coefficient of
\[
   \prod_{k\neq 0} (x \fs{F} [k]_F (\qf))
\]
is the objectionable expression $\prod_{k\neq 0} (k \qf)$; but,
as physicists say, this quantity is not directly `observable'.
For this reason, we consider the the \emph{renormalized} formal product
\[
\Theta_{F}(x;\qf) \eqdef x \prod_{k \neq 0}\frac{(x +_F [k](\qf))}
{[k](\qf)}.
\]
In section \ref{sec:prospectra} below we provide a natural setting
for such formal products.

The fixed point formula suggests that we define the equivariant
$F$-genus of $\LX$ to be
\begin{equation} \label{eq:WittenFormulaRenorm}
  \tilde{p}^{\LX}_F(1) \eqdef p^{X}_{F}\left[\prod_{j = 1}^{d}
  \frac{x_{j}}{\Theta_{F}(x_{j};\qf)} \right].
\end{equation}

\subsection{Examples} \label{sec:examples}

\subsubsection{The additive group law}

When $F$ is the additive group law $\Theta_{F}$ becomes
\begin{align*}
\Theta_{\G_a}(x,\qf)
& = x \prod_{k \neq 0} (\frac{x}{k \qf} + 1) \\
& = x \prod_{k>0}(1 - \frac{x^2}{k^2\qf^2}).
\end{align*}
This is the Weierstrass product for $\pi^{-1} \qf \sin \qf^{-1} \pi x$, so
for the theory $\CH$ of \eqref{ex:CH}, formula
(\ref{eq:WittenFormulaRenorm})
gives
\[
 \tilde{p}^{\LX}_F(1) = (2 \pi i / \qf)^{d} \left[\prod_{j = 1}^{d}
 \frac{x_{j}/2}{\sinh (x_{j}/2)}\right][X] \;.
\]
This is just the $\hat A$--genus of $X$, up to a normalization
depending on the dimension of $X$. In \cite{MFA:CircularSymmetry},
Atiyah rewrites the formal product
\[
  x \prod_{k \neq 0} (x + k \qf)
\]
as
\[
x \left(\prod_{k>0} k\right)^{2}\prod_{k > 0} \left(\frac{x^{2}}
{k^{2}} - \qf^{2}\right)
\]
and invokes zeta-function renormalization \cite{Den:det} to replace
$(\prod_{k>0}k)^{2}$ with $2 \pi$, yielding
\[
   2 \pi x \prod_{k > 0} (\frac{x^{2}}{k^{2}} - \qf^{2}) \;;
\]
specializing $\qf$ to $i$ then gives the classical expression. From our
point of view it's natural to think of the Chern class $\qf$ of $\odr{1}$
as the holomorphic one-form $z^{-1}dz$ on the complex projective line,
and thus to identify $\qf$ with its period $2 \pi i$ with respect to
the equator of $\C P_1$ as in \S 2.1 of \cite{Deligne:Galois}: the `Betti
realization' of the Tate motive $\Z(n)$ is $(2 \pi i)^n \Z \subset \C$.

\subsubsection{The multiplicative group law}
In the case of the equivariant $K$-theory $\CK$ of example
\eqref{ex:kthy}, the Euler class of a line
bundle $L$ is $\eT (L) = 1 - L$.  Writing $\qm$ for the generator
$\odr{1}$ of $\chars$, the Euler class of $\odr{1}$ is
\[
   \qf = \eT(\odr{1}) = 1 - \qm.
\]
The multiplicative group law \eqref{eq:gm} gives
\[
   \eT (L) \fs{\Gm} [k] (\qf) = 1 - \qm^{k}L,
\]
and so the formal product $\Theta_{\G_m}(x,\qf)$
becomes
\begin{equation}\label{eq:theta-gm}
\begin{split}
\Theta_{\G_m}(x,\qf)
& = (1-L)\prod_{k > 0} \frac{(1 - L \qm^{k})(1 - L \qm^{-k})}
                          {(1-\qm^{k}) (1-\qm^{-k})}  \\
& = \left(\prod_{k > 0} L\right)
    (1-L)\prod_{k > 0} \frac{(1 - L \qm^{k})(1 - L^{-1}\qm^{k})}
                       {(1-q^{k})^{2}}.
\end{split}
\end{equation}
Aside for the powers of $L$, this is essentially the product
expansion for the Weierstrass $\sigma$ function
\[
\sigma(L,q) = (1 - L) \prod_{k>0} \frac{(1 - q^k L)(1 - q^k L^{-1})}
{(1 - q^k)^2} \in \Z[L,L^{-1}][[q]]
\]
(see for example \cite{MazurTate:sigma} or p. 412 of
\cite{Silverman:EllipticCurvesII}). The infinite factor is
objectionable: in the product \eqref{eq:WittenFormulaRenorm} defining
the hypothetical $\CK$-genus of $\LX$ this factor contributes an infinite
power of $\Lambda^{\rm top} \TX$, but if $c_1 X = 0$ ({\it e.g.} if
$X$ is Calabi-Yau) and we are careful with the product, we can replace
$\Theta_{\Gm}$ with $\theta$. The resulting invariant is the Witten
genus \cite{Witten:Dir,AHS:ESWGTC}. Segal \cite{Segal:Ell} replaces
the formal product
\[
  \prod_{k \neq 0} (1 - q^{k}L) (1-q^{-k}L)
\]
which arises in the multiplicative case with
\[
  (\prod_{k>0} q^{-k} L)\prod_{k > 0} (1-q^{k}L) (1-q^{k}L^{-1}).
\]
He eliminates the infinite product of $L$'s by assuming
$\Lambda^{\rm top}\TX$ trivial and he uses zeta-function renormalization
to replace $q^{-\sum k}$ with $q^{-1/12}$.

\section{A Riemann-Roch formula for the quotient of a formal group by
a free subgroup}
\label{sec:riemann-roch}

The starting point for this paper was the discovery that the formal products
\eqref{eq:WittenFormula} and \eqref{eq:WittenFormulaRenorm} which
arise in applying the fixed point formula to study the $F$-genus of the free
loop space are \emph{precisely} the same as those obtained from the
Riemann-Roch theorem for the quotient of $F$ by a free cyclic
subgroup.  We explain this in section \ref{sec:RiemannRochApplied},
after briefly reviewing finite quotients of formal groups, following
\cite{Lubin:Fin,Ando:PowerOps}.

\subsection{The quotient of a formal group by finite subgroup}
\label{sec:lubin}

In this section, we assume that $F$ is a formal group law over a complete
local domain $R$ of characteristic $0$ and residue characteristic
$p>0$. If $A$ is a complete local $R$-algebra, the group law
$F$ defines a new abelian group structure on the maximal ideal $\m{A}$
of $A$. We will refer to $(\m{A},+_F)$ as the group $F(A)$ of $A$-valued
points of $F$.

If $H$ is a \emph{finite} subgroup of $F(R)$, then Lubin shows that
there is a formal group law $F/H$ over $R$, determined by the
requirement that the power series
\begin{equation}\label{eq:lubin}
     f_{H} (x) \eqdef \prod_{h\in H} (x\fs{F}h) \in R\psb{x}
\end{equation}
is a homomorphism of group laws $F\xrightarrow{f_{H}}F/H$; in other
words there is an equation
\[
     F/H (f_{H} (x),f_{H} (y)) = f_{H}\left(F (x,y)\right).
\]
The main point is that the power series $f_H$ is constructed so
the kernel of $f_H$ applied to $F(R)$ is the subgroup $H$.

The coefficient $f_{H}' (0) = \prod_{h\neq 0}h$ of the linear term of
$f_{H} (x)$ is not a unit of $R$, and so $f_{H}$ is an isomorphism of
formal group laws only over $R[f'_H(0)^{-1}].$  Over this ring,
we might as well replace $f_H$ with the strict isomorphism
\[
     g_{H} (x) \eqdef x \prod_{h \neq 0} \frac{x\fs{F}h}{h} =
\frac{f_{H} (x)}{f_{H}' (0)} ,
\]
and define $G$ to be the formal group law
\[
   G (x,y) = g_{H}\left(F (g_{H}^{-1} (x),g_{H}^{-1} (y))\right)
\]
over $R[f'_H(0)^{-1}]$; then $F/H$ and $G$ are
related by the isomorphism
\[
    t (x) = f'_{H} (0)x.
\]

If $F$ is the group law and $R$ the ring of coefficients of
a parametrized theory, then the Riemann-Roch formula
\eqref{eq:RiemannRoch} for a compact complex-oriented manifold $X$
is the equation
\begin{equation}\label{eq:rr-subgp}
      p_{G}^{X} (u) =
      p_{F}^{X} \left[u\prod_{j=1}^{d} \frac{x_j}{g_H (x_j)}\right]
\end{equation}
over $R[f'_H(0)^{-1}]$.

\subsection{The case of a free cyclic subgroup}
\label{sec:RiemannRochApplied}

Now suppose that $(E,x,F)$ is a parametrized theory.  Let $R$ be the
even periodic theory defined by the formula
\[
  R (X) \eqdef E (B\T\times X).
\]
The projection $B\T\to \point$ gives a natural transformation $E (X)\to
R (X)$.  In particular the coordinate $x\in E (P)$ gives a coordinate
$x\in R (P)$.  The group law $\bF$ is just the group law
$F$, considered over the $E$-algebra $R$.

What extra structure is available over $R$?  A
character $\lambda\in \chars$ gives a map $B\T \to P$, and so an
$R$-valued point $u(\lambda)$ of $\gpof{E}$.  As $\lambda$ varies
through the group of characters, these points assemble into a
homomorphism of groups
\begin{equation}\label{eq:univ-hom-chars-to-F}
  u\: \chars \to \bF(R).
\end{equation}
It is easy to see that this is the inclusion of a subgroup.  As usual
we write $\qf=u (1)$ for the generator of $\chars$ using the
isomorphism $\chars\cong\Z$.  We write $(\qf)$ for the subgroup $u
(\chars)$ generated by $\qf$.

The analogue of Lubin's formula \eqref{eq:lubin} is
\begin{equation} \label{eq:fq}
  f_{(\qf)}(x) = \prod_{k\in \Z} (x\fs{F} [k](\qf)),
\end{equation}
and the resulting Riemann-Roch formula
\[
  \pr{X}_{\bF/(\qf)} (u) =
  \pr{X}_{F}\left[ \;u\; \prod_{j=1}^{d} \prod_{0 \neq k}
    \frac{1}{x_{j} \fs{F} [k]_{F}(\qf)}\right]
\]
has right-hand side identical to the formal product
\eqref{eq:WittenFormula} arising from the fixed-point formula.

Of course this does not avoid the problems of applying the
fixed-point formula to the free loop space. From this point of view,
the trouble is that the quotient object $\bF/(\qf)$ isn't a formal group;
for example the coefficient of $x$ is the product \eqref{eq:fq} is
$\prod_{k\neq 0}(k \qf)$  This particular problem is fixed by using
\begin{align*}
    g_{(\qf)} (x)
& = x \prod_{0\neq k\in \chars}
    \frac{x \fs{F} [k]_{F} (\qf )}{[k]_{F} (\qf)} \\
& = \Theta_{F}(x,\qf).
\end{align*}
The Riemann-Roch formula in this case is
\begin{equation}\label{eq:RR-g}
\pr{X}_{G} (u) = \pr{X}_{F}\left[ u \prod_{j=1}^{d}
    \frac{x_j}{\Theta_{F}(x_j,\qf)}\right]
\end{equation}
with right-hand-side identical to the renormalized genus
\eqref{eq:WittenFormulaRenorm}.

However, expanding such Weierstrass products as formal power
series is still highly nontrivial, as the examples above have shown.

\section{The structure of $\bF/ (\qf)$}
\label{sec:GeneralizedTateCurve}

Suppose now that $F$ is a formal group law over a complete local
ring $E$ with residue field $k$ of characteristic $p > 0$, and
that $\bF$ is its pullback over the power series algebra $R=E\psb{\qf}$.
The coordinate defines a homomorphism
\begin{equation}\label{eq:univ-hom-Z-to-F}
\begin{split}
 u \: \Z & \to \bF (R) \\
      u(n) & = [n]_{\bF} (\qf);
\end{split}
\end{equation}
but there is no reason to expect the cokernel of this homomorphism to
be a formal group. However, $\bF/(\qf)$ certainly makes sense
as a group-valued functor of complete local $R$-algebras.

In section \ref{sec:katz-mazur}, we show by construction that the
torsion subgroup of $\bF/ (\qf)$ has a natural approximation by a
representable functor.  We construct a formal group scheme
$\Tate{F}$ over $R$ together with a natural transformation
\[
     \Tate{F}\torsion \to \bF/(\qf)\torsion
\]
of group-valued functors, which is an isomorphism if $\qf$
has infinite order in $F$. The formal group scheme $\Tate{F}\torsion$
is our model for $\bF/(\qf)$.

Because $E$ has finite residue characteristic, we work $p$-locally:
the formal group scheme $\Tate{F}\ptorsion$ is a $p$-divisible group
in the sense of  \cite{Tate:pdivisiblegps,Demazure:pdivisiblegps}. The
$p$-torsion subgroup
$G\ptorsion$ of a formal group $G$ of finite height is a
\emph{connected} $p$-divisible group, but $\Tate{F}\ptorsion$ is not;
its maximal \'etale quotient is a constant height-one $p$-divisible group
$\Qp/\Zp$, and its connected component is just the $p$-divisible group
$F\ptorsion$ of $F$. In other words, there is an extension
\[
    \bF\ptorsion \xra{} \Tate{F}\ptorsion \xra{} \Qp/\Zp
\]
of $p$-divisible groups over $R$. We will see in \eqref{t-le:univ-ext}
that it is in fact the \emph{universal} example of such an extension, and
it follows that if $E_{n}$ is the ring which classifies lifts of a
formal group $G$ of height $n$ over an algebraically closed field $k$
\cite{LubinTate:FormalModuli}, then $R = E_{n}\psb{\qf}$ represents the
functor which classifies lifts of a $p$-divisible group of height $n+1$
with connected component $G\ptorsion$.   Thus $\qf$ may be viewed as a
``Serre-Tate parameter'' in the sense of \cite{Katz:SerreTate}.

\subsection{A model for the torsion subgroup of $F/ (\qf)$}
\label{sec:katz-mazur}

One difficulty in representing the quotient $\bF/ (\qf)$ over $E\psb{\qf}$ is
that the subgroup $(\qf)$ does not act universally freely on
$\bF$: consider, for example, the specialization $\qf=0$.
It turns out that, as long as one restricts to the torsion subgroup of
$\bF/ (\qf)$, this is the only obstruction.  By freeing up the
action, we are able to construct a representable functor whose torsion
subgroup coincides with that of $F/ (\qf)$ whenever $\qf$ is of
infinite order.

Following \cite[\S 8.7]{KatzMazur:Ari}, let $\Tate{F}$ be the scheme over
$R=E\psb{\qf}$ defined by the disjoint union
\[
\Tate{F} = \bigcup_{a \in \Q \cap [0,1)}
F_R \times \{a\}.
\]
If $A$ is a complete local $R$-algebra, then
\[
\Tate{F}(A) =  \{ \text{pairs }(g,a)\text{ with }
g \in F(A)\text{ and } a \in \Q\cap [0,1). \}
\]
This has a group structure given by
\[
(g,a)\cdot (h,b) = \begin{cases}
     (g \fs{F} h,a + b) & \text{ if $a + b < 1$} \\
     (g \fs{F} h \fm{F}\qf, a + b - 1) &
                            \text{ if $a + b\geq 1$.}
                                      \end{cases}
\]
By construction,  $\Tate{F}(A)$ is the quotient in the exact sequence
\begin{equation} \label{eq:Texactsequence}
\begin{array}{clrl}
 0 \xra{} & \Z \xra{} & F(A) \times  \Q
                      & \xra{} \Tate{F}(A)\xra{} 0  \\
                    & n \mapsto  & ([n](\qf),n)  \\
                    &            & (x,a) &
              \mapsto(x \fm{F} [\flat(a)](\qf),\sharp(a)),
\end{array}
\end{equation}
where $\flat(a)$ and $\sharp(a)$ are the integral and
fractional parts of the rational number $a$.  Equivalently,
$\Tate{F}$ is the pushout of $\Z \to \Q\to
\Q/\Z$  along the homomorphism $u$:
\[
\begin{CD}
\Z @>>> \Q @>>> \Q/\Z\\
@V u VV  @VVV          @VV = V \\
F_{R} @>>> \Tate{F} @>>> \Q/\Z.
\end{CD}
\]
Thus $\Tate{F}$ is a kind of homotopy quotient of $\bF$ by $\Z$.
\begin{Proposition} \label{Th:TeqFmodq}
Projection onto the first factor in the construction above defines
a natural transformation
\[
   \Tate{F}\torsion \to (\bF/ (\qf))\torsion,
\]
which is an isomorphism if $\qf$ is not in $F (A)\torsion$.
\end{Proposition}

\begin{proof}
Let $A$ be a complete local $R$ algebra, and suppose that $\qf$ is not
torsion in $F (A)$.  We see from
(\ref{eq:Texactsequence}) that there is an isomorphism
\[
\Tate{F} (A)\cong \bF/ (\qf) (A) \times \Q
\]
which is compatible with the projection and clearly induces an
isomorphism
\[
\Tate{F}\torsion (A)\cong (\bF/ (\qf))\torsion (A).
\]
\end{proof}

\begin{Remark}
The Proposition holds in any context in which $\qf$ is not
torsion.  Another source of examples is $p$-adic fields.  Let
$S=E[\tfrac{1}{p}]\lsb{\qf}$.  Let $L$ be a complete
nonarchimedean field with norm $\padnrm{\slot}$.  Given a continuous
map $S\to L$, the group law $\bF$ defines a group structure on the set
$\{v\in L| \padnrm{v}<1\}$; and we write $\bF (L)$ for this group.
The subgroup generated by $\qf$ is necessarily free, and
the argument shows that there is an isomorphism of groups
\[
    \Tate{F}\torsion (L)\cong (\bF/ (\qf))\torsion (L).
\]
\end{Remark}

\begin{Remark}
If the ring $E$ is the ring of coefficients in an even periodic ring theory,
then the preceding construction could be carried out with $R=E (B\T)$
and the homomorphism $u$ \eqref{eq:univ-hom-chars-to-F}.  The result is an
extension
\[
\gpof{E} \rightarrow \Tate{\gpof{E}} \rightarrow (\chars\otimes \Q)/\chars
\]
of group schemes over $R$.  Proposition \ref{Th:TeqFmodq} becomes an
isomorphism
\[
     \Tate{\gpof{E}}\torsion \cong ((\gpof{E})_{R}/\chars)\torsion.
\]
\end{Remark}

\subsection{Notation}
\label{sec:GG}
If $A$ is  an abelian group, then $A_{T} =  \spec T^{A}$ is the
resulting constant formal group scheme over $T$.
The category of \emph{test rings} is the category of Artin local
$E$-algebras with residue field $k$.

\begin{Definition} \label{def:hom}
If $G$ is a formal group scheme over a ring $T$, and $A$ is an abelian
group, let $\Hom[A,G]$ be the group of homomorphisms
\[
     A_{T} \xra{} G
\]
of formal group schemes over $T$.  Similarly, let $\Ext[A,G]$
be the set of isomorphisms classes of extensions of formal group schemes
\[
     G \xra{} X \xra{} A_{T}.
\]
If $G$ is a formal group scheme over $E$, let $\uHom[A,G]$ and
$\uExt[A,G]$ be the functors from test rings to groups such that
\[
     \uHom[A,G] (T) \eqdef  \Hom[A_{T},G_{T}]
\]
and
\[
     \uExt[A,G] (T) \eqdef \Ext[A_{T},G_{T}].
\]
\end{Definition}

Now if $G$ is a formal group over $E$, pulling back over $G$
defines a natural point
\[
       \Delta\:  G \xra{} G\times G = G_{G}
\]
and hence a homomorphism
\[
       u\: \Z_{G} \xra{} G_{G}.
\]
It is clear that this gives an isomorphism $G\cong\uHom[\Z,G]$.

Equivalently, if $F$ is a formal group law over $E$, then
$R=E\psb{\qf}$ pro-represents the functor $\uHom[\Z,F]$ on the
category of test rings, with universal example $u (n) = [n]_{F} (\qf)$
\eqref{eq:univ-hom-Z-to-F}.  Similarly, if $E$ is the ring of
coefficients  of an even periodic ring theory, then $R = E(B\T)$
pro-represents the functor $\uHom[\chars,\gpof{E}]$, with the
homomorphism $u$ of \eqref{eq:univ-hom-chars-to-F} as the universal example.

\subsection{Universal properties}
\label{sec:universalprops}
\subsubsection{A universal extension}

A  continuous homomorphism of $E$-algebras from $R$ to a test algebra
$T$ defines an extension
\[
     F_{T} \xra{} (\Tate{F})_{T} \xra{} (\Q/\Z)_{T}
\]
and hence an extension
\[
     F\ptorsion_{T} \xra{} (\Tate{F}\ptorsion)_{T} \xra{} (\ZpmZ)_{T} \;.
\]
of torsion subgroups.

\begin{Lemma} \label{t-le:univ-ext}
The ring $R$ pro-represents the functor $\uExt[\Qp/\Zp,F]$, with
$\Tate{F}$ as the universal example.
\end{Lemma}

\begin{proof}
In the exact sequence
\[
             \uHom[\Q,F]  \xra{}
             \uHom[\Z,F] \xra{}
             \uExt[\Q/\Z,F] \xra{}
             \uExt[\Q,F],
\]
the first and last terms  are zero because $p$ acts
nilpotently on $F$ and as an isomorphism on
$\Q$.
\end{proof}

\subsubsection{A universal $p$-divisible group}

Any $p$-divisible group $\Gamma$ over a field $k$ is naturally an extension
\begin{equation}\label{eq:comp-exact-seq}
   \Gamma^0 \xra{} \Gamma \xra{} \Gamma_{\etale}
\end{equation}
of a connected group by an \'etale group. If the residue field $k$
is algebraically closed, then the sequence
\eqref{eq:comp-exact-seq} has a canonical splitting
\cite[p. 34]{Demazure:pdivisiblegps}. We will be interested in
the case when $\Gamma_{\etale}$ has height one, which is to say
that it is isomorphic to the constant group scheme $(\Qp/\Zp)_k$.

Tate showed \cite{Tate:pdivisiblegps} that the functor $G\mapsto
G[p^{\infty}]$ is an equivalence between the categories of formal
groups of finite height and connected $p$-divisible groups. Let's fix
a one-dimensional formal group $\Gamma^{0}$ of height $n$ over the
algebraically closed field $k$ and define $\Gamma$ to be the product
extension
\[
    \Gamma^0 \to \Gamma \to  (\T^* \otimes \Qp/\Zp)_k \;.
\]

\begin{Definition} \label{def:lifts}
If $G$ is a $p$-divisible group over $k$, and if $T$ is a test ring,
then a \emph{lift} of $G$ to $R$ is a pair $(F,\delta)$
consisting of a $p$-divisible group $F$ over $R$ and an isomorphism
\[
     F_{k} \xrightarrow[\cong]{\delta} G
\]
of $p$-divisible groups over $k$.  An equivalence of lifts
$(F,\delta)$ and $(F',\delta')$ is an
isomorphism $f\: F\to F'$ such that
\[
       \delta = \delta' f_{k}.
\]
The set of isomorphism classes of lifts of $G$ to $R$ will be denoted
$\LT{G} (R)$.  As $R$ varies, $\LT{G} (R)$ defines a functor from test
rings to sets.
\end{Definition}

Lubin and Tate construct a formal power series algebra $E_n$ over the Witt
ring of $k$ which pro-represents the functor $\LT{\GammaF}$.  There is
an even periodic cohomology theory with $E_n$ as ring of coefficients,
and the universal lift $F$ of $\Gamma^0$ as formal group.

\begin{Theorem} \label{t-th:univ-def}
The ring $R=E_n(B\T)$ pro-represents $\LT{\Gamma}$, with universal example
$\Tate{F}\ptorsion$.
\end{Theorem}

\begin{proof}
Let $T$ be a test ring.  Suppose that $(H,\epsilon)$ is a lift of
$\Gamma$.  Then $(H^{0},\epsilon^{0})$ is a lift of $\Gamma^{0}$.
According to \cite{LubinTate:FormalModuli}, there is a unique pair
$(f,a)$ consisting of a map $f\: E_n \to T$ and an isomorphism
\[
    a\: (F_{T},\epsuniv_{T})\cong (H^{0},\epsilon^{0})
\]
of lifts of $\GammaF$. On the other hand, $\epsilon_{\etale}$ induces
an isomorphism
\[
     H_{\etale} = (H_{\etale} (k))_{T}\cong
     (\Gamma_{\etale} (k))_{T} = (\chars \otimes \Qp/\Zp)_{T}.
\]
Assembling these gives an extension
\[
       F_{T}\rightarrow H \rightarrow (\chars \otimes \Qp/\Zp)_{T}
\]
which defines an isomorphism
\[
      \LT{\Gamma} \cong \uExt[F,\chars \otimes \Qp/\Zp] \;.
\]
$R$ pro-represents the right-hand side by Lemma \ref{t-le:univ-ext}.
\end{proof}

\begin{Remark} \label{rem:KatzMazur}
The analogous result for the ordinary multiplicative group
$\rGm$ is described in \cite[\S 8.8]{KatzMazur:Ari}. Closely
related examples occur in \cite{AMS:Tate}, which are motivated
by purely homotopy-theoretic questions about Mahowald's root
invariant.\end{Remark}

\subsection{Examples}

If $L \in A((q))$ is a unit, the formal product defining the Weierstrass
function of\ref {sec:examples} defines an element $\sigma(L,q) \in A((q))$.
The functional equation $\sigma(qL,q) = (-L)^{-1} \sigma(L,q)$ then
implies that the modified product
\[
\sigma [L,r] \; \eqdef \; q^{-\flat(r)(\flat(r) + 1)/2}
(-L)^{\flat(r)} \sigma(L,q)
\]
(where $\flat(r)$ is the integral part of $r \in \Q$) satisfies
the identity $\sigma[qL,r+1] = \sigma[L,r]$ and can thus be regarded
as a function from a localization of $\Tate{\G_m}$ to the $\Z((q))$-line.
Since $\sigma(L,0) = 1 - L$ we can think of the modified function
as a deformation of the usual coordinate at the identity of the
multiplicative group.

Similarly, if we regard $\qf$ as an element
of the locally compact field $\C$, the modified sine function
\[
s[x,r] \; \eqdef \; \pi^{-1} \qf \sin \qf^{-1} (\pi x - r)
\]
satisfies the identity $s[x + \qf,r+1] = s[x,r]$ and so can be interpreted
as a function from $\Tate{\G_a}(A)$ to $A$. It is also a deformation of
the usual coordinate on the additive group, in that $s[x,0] \to 0$
as $\qf \to \infty$.

In general, however, there seems to be no reason to expect that a
coordinate on $F$ will extend to a coordinate on $\Tate{F}$: our
construction yields a group object, but not a group \emph{law}. In
the two examples above, we do have (something like) coordinates,
which define interesting genera: ordinary cohomology leads to the
$\hat A$-genus, suitably normalized, if $r=0$; but if $r \neq 0 \in \Q/\Z$,
Cauchy's theorem (applied to a small circle $\bf C$ around the origin)
yields
\[
\frac{1}{2\pi i} \int_{\bf C} p^X_{\G_a} \left[ \prod_{j=1}^{j=d} \frac
{\qf^{-1} \pi x_j z}{\sin \qf^{-1} \pi (x_j z + r \qf)} \right] z^{-d-1} \; dz
= (\frac{\pi}{\qf \sin \pi r})^{d} \chi(X) \;.
\]
The $K$-theory genus extends similarly, but the resulting function is
just $(-1)^d q^{-d \flat(r)(\flat(r) + 1)/2}$ times the standard elliptic
genus.

Any multiplicative cohomology theory $E$ can be described as taking
values in a category of sheaves over $\spec E(\point)$, and the
Borel extension of such a theory takes values in sheaves over
$\spec (E(\point) \psb {\qf})$. The construction in \ref{sec:katz-mazur}
of the Tate group as a disjoint union of copies of such affines
implies that a theory $E$ with formal group $F$ has a natural
extension to an equivariant theory taking values in a category
of sheaves over the group object $\Tate{F}$. Similarly, a suitable
localization of the Borel extension defines an equivariant theory
taking values in sheaves over a suitable localization of $\Tate{F}$.
This resembles (but is easier than) the constructions of
\cite{GKV:ell,Grojnowski:ell,KnutsonRosu}, for here we're only
patching together Borel extensions.]

\section{Prospectra and equivariant Thom complexes}
\label{sec:prospectra}

Cohen, Jones, and Segal show that pro-objects in the
category of spectra are the appropriate context in which to study Thom
complexes, and so umkehr maps, for semi-infinite vector bundles. In
this section we observe that the ideas of this paper fit very naturally
into their framework.

\subsection{Thom prospectra}

If $V$ and $W$ are complex vector bundles over a space $X$, we extend
the notation for Thom isomorphisms in \ref{rem:thom-isos} by writing
\[
\alpha^{V}_{W}  \eqdef \alpha^{V} \circ (\alpha^{W})^{-1} \:
                   E(X^{W}) \rightarrow E(X^{V}) \;.
\]
`Desuspending' by $V \oplus W$ then gives a homomorphism
\[
 \alpha^{-W}_{-V} \: E(X^{-V}) \rightarrow E(X^{-W}) \;.
\]
The inclusions of a filtered vector bundle
\[
 {\V} : 0 = V_{0} \subset V_{1} \subset \cdots
\]
define maps
\[
i_{V_{n}} : X^{V_{n-1}} \rightarrow X^{V_{n}}
\]
of Thom spectra, which desuspend to define a pro-object
\[
 X^{-{\V}}  \eqdef  \{ \dots \rightarrow X^{-V_{n}} \xrightarrow
{i^{V_n}} X^{-V_{n-1}} \rightarrow \dots \}
\]
in the category of spectra. The $E$-cohomology of $X^{-{\V}}$ is the colimit
\[
E(X^{-{\V}})\eqdef \colim_{n} E (X^{-V_{n}}),
\]
as in the appendix to \cite{CJS:Floer}.

\begin{Lemma} \label{t-le:E-of-i-n}
On cohomology, the homomorphism induced by $i^{V_{n}}$ is $e(V_{n}/V_{n-1})
\alpha_{-V_{n-1}}^{-V_{n}}$.
\qed
\end{Lemma}

\begin{Example}
For any integer $n$, let
\[
I_{n}  =  \bigoplus_{n \geq |k| >0} \odr{k};
\]
thus $\IV{} = \colim_{n} I_{n}$ is a filtered $\T$-vector bundle over
a point. More generally, if $V$ is a complex vector bundle over a space
$X$, let
\[
\IV{V}  \eqdef  V\otimes \IV{}
\]
be the corresponding filtered $\T$-vector bundle. In this notation
the Laurent approximation \eqref{eq:normalbundle} to the formal normal
bundle of the constant loops in $\LX$ is $\IV{\TX}$.

We write $j_{n}$ for the map
\[
    j_{n} = i^{\In{n}{V}} \: X^{-\In{n}{V}} \to X^{-\In{n-1}{V}} \;.
\]
\end{Example}

\begin{Example}
If $V_n$ is the sum of $n$ copies of $\odr{1}$, considered as a $\T$-vector
bundle over a point, then the Borel construction on ${\V}$ is the Thom
prospectrum for $\CP_{\infty}$ constructed in \cite{CJS:Floer}.
\end{Example}

\subsection{Equivariant cohomology}

Suppose now that $\hET$ is a suitable extension of an equivariantly oriented
theory $\ET$, and let $X$ be a finite complex with trivial $\T$-action, as
above. As $n$ varies, the Thom isomorphisms
\[
      \alpha^{\In{n}{V}} \: \ET(X) \xrightarrow[\cong] \ET(X^{-\In{n}{V}})
\]
are not compatible with the maps $j_{n}$; but this can be cured over $\hET$
by a suitable renormalization. By \eqref{t-le:localization} the class
\[
u_{n} (V) \eqdef e (V \otimes I_{n}) = \prod_{0<|k|\leq n}
\eT (V\otimes \C (k)).
\]
is a unit of $\hET (X)$, and the homomorphism
\[
\omega_{n} (V)\: \hET (X) \xrightarrow{} \hET (X^{-\In{n}{V}})
\]
defined by the formula $\omega_{n} (V) =  u_{n} (V)\alpha^{-\In{n}{V}}$
is an isomorphism.
\begin{Theorem} \label{t-th:ind-thom}
The diagram
\[
\begin{CD}
  \hET (X) @= \hET (X) \\
  @V {\omega_{n-1} (V)} VV @VV {\omega_{n} (V)} V \\
  \hET (X^{-\In{n-1}{V}}) @> j_{n} >> \hET (X^{-\In{n}{V}})
\end{CD}
\]
commutes; in particular, the maps $\omega_{n} (V)$ assemble into a
``Thom isomorphism''
\[
      \omega^{V} \: \hET (X) \xra{\cong} \hET (X^{-\IV{V}}) \;.
\]\qed
\end{Theorem}

If $V$ and $W$ are two complex vector bundles over $X$, then
(as with the usual Thom isomorphism) we define
\[
   \omega_{W}^{V} \: \hET (X^{-\IV{W}}) \xra{\cong} \hET (X^{-\IV{V}}).
\]
\begin{Corollary} 
If the vector bundle $V$ has rank $d$, then the relative isomorphism
\[
     \omega^{V}_{d}\: \hET (X^{-\IV{\triv{d}}}) \rightarrow \hET (X^{-\IV{V}})
\]
is given by the formula
\[
      \omega^{V}_{d} =
     \left(\prod_{j=1}^{d}\prod_{k\neq 0} \frac{x_{j} \fs{F}[k] (\qf)}
                                     {[k] (\qf)} \right)
                               \alpha^{-\IV{V}}_{-\IV{\triv{d}}},
\]
where the $x_{j}$ are the terms in the formal factorization
\eqref{eq:chern-splitting} of the total $F$-Chern class of $V$.\qed
\end{Corollary}

The terms in this product are well-defined at any finite stage, if
not in the limit. If $\ET$ is $\CH$ or $\CK$ (and $c_1(V) = 0$), then
the infinite products converge and we have the
\begin{Corollary} \label{t-co:procomplex-coh-thom}
If the vector bundle $V$ has rank $d$, then the diagram
\[
\begin{CD}
    \hET(X^{\triv{d}}) @> {\alpha^{-\IV{\triv{d}}}_{\triv{d}}} >>
    \hET(X^{-\IV{\triv{d}}}) \\
    @V{\hat{\alpha}^{V}_{\triv{d}}} VV @VV {\omega_{\triv{d}}^{V}} V \\
    \hET(X^V) @> {\alpha^{-\IV{V}}_V} >> \hET(X^{-\IV{V}}) \;.
\end{CD}
 \]
commutes; where
\[
\hat{\alpha}^{V}_{\triv{d}} \; \eqdef \; \frac{\Theta_F(V)}{c_d(V)} \;
\alpha^V_{\triv{d}} \;.
\]
\end{Corollary}

More precisely, the assertion is that these examples are naturally
topologized, and that at any finite stage the diagram commutes modulo
error terms which converge in this topology to zero as $n$ grows.
\qed

In other words, in a suitable theory with group law $F$ the isomorphism
$\omega^{V}_{\triv{d}}$ is very much like a Thom isomorphism
$\hat{\alpha}^V_{\triv{d}}$ for a theory with $F$ replaced by its
extension $\Tate{F}$. From this point of view, the fact that Witten's
formula for the $\alpha$-genus of $\LX$ equals the $\hat \alpha$-genus
of $X$ can be interpreted as saying that the inclusion $j \: X \to \LX$
behaves as if there is a cohomological analogue
\[
 j_* \: \hET(X^{-\IV{TX}}) \to \hET(\LX)
\]
of the Thom collapse map, with an associated umkehr $j_F = j_* \omega^{TX}$
satisfying the projection formula
\[
  j_F j^*(x) = x \cdot \eT(\nu) \;.
\]
Perhaps the intuition underlying the physicists' interest in
elliptic cohomology is that a reasonable equivariant theory applied
to the geometrical object $X^{-\IV{TX}}$ captures more information
than that theory does, when applied directly to $X$. Thus the equivariant
$K$-theory of this formal neighborhood of $X$ is the Tate elliptic
cohomology of $X$. This seems to be related to the recent construction
(by Kontsevich and others) of new invariants for singular complex
algebraic varieties, by considering them as varieties over the
Laurent series field $\C((q))$.

\subsection{Polarizations}

In the preceding account, the role of the rational parameter
in the construction of $\Tate{F}$ is geometrically unmotivated, because we
have ignored some issues connected with the polarization
\cite{CJS:Floer} of the loopspace.

Such a structure is an equivalence class of splittings of the tangent bundle
of $\LX$ into a sum of positive- and negative-frequency components: if
$X$ is complex-oriented ({\it e.g.} symplectic, with a choice of compatible
almost-complex structure), then the composition
\[
\LX \rightarrow {\mathcal L}(BU) = U \times BU \rightarrow U/SO
\]
of the map induced by the classifying map for the tangent bundle of $X$ with
the projection to a classifying space for such splittings defines a canonical
polarization.

\subsubsection{The universal cover of the free loopspace}

If $X$ is simply-connected, the fundamental group of its free loopspace
will be isomorphic to $H_2(X,\Z)$, which will be nontrivial in general.
There is thus good reason to consider the simply-connected cover $\OLX$
of the loopspace: this can be defined as the space of smooth maps
of a two-disk to $X$, modulo the relation which identifies two maps if their
restrictions to the boundary circle agree, and if furthermore their
difference, regarded as an element of the deck-transformation group
$\pi_2(X) = H_2(X,\Z)$, is null-homotopic. The
circle acts on $\OLX$ by rotating loops, as does the fundamental
group of $\LX$, and in general the fixed-point set
\[
\OLX^{\T} \cong H_2(X,\Z) \times X
\]
will have many components. The choice of a basepoint defines a lift of
the canonical polarization to a map
\[
\OLX \rightarrow \R \times SU/SO
\]
which restricts to a locally constant map
\[
H_2(X,\Z) \times X \rightarrow \R \times SU/SO \rightarrow \R
\]
defined by evaluating $\alpha \in H_2(X,\Z)$ on the first Chern
class of $X$. We think of the polarization as defining the zero-frequency
modes in the Fourier decomposition of small loops near a fixed-point
component, so that shifting by $\beta \in H_2(X,\Z)$ gives an isomorphism
\[
T\OLX^{\T}_{x,\alpha} \cong T\OLX^{\T}_{x,\alpha + \beta}
\otimes \odr{\langle c_1(V),\beta \rangle} \;.
\]

\subsubsection{$\T$-equivariant Picard groups and orientations}

We can regard the cohomology of a space $Y$ as a sheaf of
rings over the zero-dimensional scheme
\[
  (\pi_0 Y)_{\Z} \eqdef \spec H^0(Y,\Z)  \;.
\]
Similarly, the set of equivalence classes of $\T$-line
bundles over a $\T$-space $X$ is naturally isomorphic to
$H^2(Y_{\T},\Z)$, which can be interpreted as a sheaf of
groups over $(\pi_0 Y)_{\Z}$, the fiber above component
$Y_i$ being the constant group scheme
\[
{\rm coker} \;[H^2_{\T}(\point,\Z) \to H^2_{\T}(Y_i)] \;;
\]
when the circle action on $X$ is trivial, this is just a complicated
way of indexing the summands of $H^2(Y,\Z)$.

An orientation on a complex-oriented cohomology theory $E$
with formal group $F$ defines a natural homomorphism from the
Picard group of line bundles over $Y$ to the group of $E(Y)$-valued
points of $F$. This suggests that a complex orientation for
the Tate extension of an equivariantly complex-oriented theory
should be defined as a natural transformation from the Picard
group of $\T$-line bundles over $X$ to $\Tate{F}(E(Y))$, both
regarded as schemes over $(\pi_0 Y)_{\Z}$; by restriction such a
transformation defines a map
\[
  \pi_0 Y \rightarrow \Q/\Z \;.
\]
In other words, such a generalized orientation assigns to
a component of $Y$ and a $\T$-line bundle over it, a characteristic
class in the cohomology of the component, together with an $r \in \Q$
depending on the component, which shifts by an integer when
the bundle is twisted by a character.

When $Y$ is $\OLX^{\T}$, the polarization defines a very natural map
of this sort, which sends $\alpha$ to $\langle c_1(V),\alpha \rangle$.
In more general situations, {\it e.g.} in Givental's work
\cite{Voisin:Givental} on the quantum cohomology of a symplectic manifold
$(X,\omega)$ the evaluation map $\alpha \mapsto \langle \omega, \alpha
\rangle$ plays a similar role; but in this generalization, $\omega$ no
longer needs to be an integral class.

\section{Concluding remarks}

The parts of this paper which deal principally with the fixed point
formula on the free loop space are formulated in terms of a
general equivariant cohomology theory $\ET$, but those which relate
to the quotient $\bF/ (\qf)$  use only the group law $\bF$ of the
theory $R (X)=E (B\T\times X)$, and so essentially use the Borel
theory $\borel{E}$.  We do this because in that case we can
be more specific about our constructions.  In good cases one can hope
to do better.

Specifically, suppose that the multiplication homomorphism
$\T\times\T\to \T$ induces a map
\[
    \ET \to E_{\T\times\T}
\]
and that one has an isomorphism $E_{\T\times\T} (\point)\cong\ET
(\point)\tensor{E} \ET (\point)$.  The upshot is a group scheme
$G=\spec \ET (\point)$ over $E$, such that the formal group $\gpof{E}$
associated to $E (P)$ is the completion of $G$.
Indeed \cite{GKV:ell,Grojnowski:ell}, one expects in general that there
is an abelian group scheme $G$ which is the more fundamental object in
$\T$-equivariant $E$-theory, with $\ET(\point)$ as structure sheaf.

In any case, as remarked in \S \ref{sec:GG}, over $G_{G} = G\times G$
there is a natural map
\begin{equation}\label{eq:u-again}
  \Z_{G} \to G_{G},
\end{equation}
and the natural map $G_{G} \to \Hom[\Z,G]$ is an isomorphism.  One
could consider the group $\Tate{G}$ over $G$, fitting into a
short exact sequence
\[
   G_{G} \xra{} \Tate{G} \xra{} \Q/Z.
\]
The group $\Tate{\gpof{E}}$ considered in the main text would then
arise from $\Tate{G}$ by completing in both copies of $G$.

This is the situation in $K$-theory.  Over the ring $\KT=\Z[q,q^{-1}]$ one has a homomorphism
\begin{align*}
   u\: \Z& \to{\rGm}_{\KT}\\
    n&\mapsto q^{n}
\end{align*}
and this gives rise to a group $\Tate{\rGm}$ over $\KT$.  Over the
completion $\KT\to K (B\T) = \Z\psb{\qf}$ one recovers the group
$\Tate{\Gm}$ considered in the main text.  On the other hand, the
group $\Tate{\rGm}\torsion$ becomes isomorphic to the torsion subgroup
of the classical Tate curve $\Tate{q}$
already over the suitable localization $\CK=\Z \lsb{q}$ of
\eqref{ex:kthy}, where $u$ is the inclusion of a sub-groupscheme.  See
\cite[\S 8.8]{KatzMazur:Ari} for details.

Interestingly enough, in the case of $K$-theory this \emph{is} the
solution to the problem of the infinite products.  We have
purposefully written the equivariant Euler class of $L\otimes \odr{k}$ as
$e (L) \fs{F} [k]_{F} (\qf)$, as in
any finite situation the formula
\[
    1 - L_{1} L_{2} = (1 - L_{1}) + (1-L_{2}) - (1-L_{1}) (1-L_{2})
\]
makes it possible to calculate the Euler class of a vector bundle
using the formal multiplicative group, i.e. the right-hand side.
However, in order to calculate the infinite product $\Theta_{\Gm}
(x,\qf)$ \eqref{eq:theta-gm}, one is \emph{forced} to use the
left-hand side.  Thus from this point of view the renormalization is
handled by knowing about the global multiplicative group, instead of
merely its formal completion.

It seems reasonable to hope that the rich theory of $\T$-equivariant
ring spectra will provide additional examples of
such $\ET$.  If $\ET$ is a complex-oriented $\T$-equivariant ring
spectrum, then the cohomology $\ET (\C \P)$ (where $\C \P$ is the
space of lines in the ambient $\T$-universe) carries the structure of
a $\T$-equivariant formal group law (in the sense of
\cite{CGK:EFG}).  This is, among other things, a formal group law $F$ and
a homomorphism
\[
  v\: \chars \to F (\ET (\C\P)).
\]
In the situation we are considering, $F$ is the completion of the
group $G=\spec \ET$, and the homomorphism $v$ is obtained from the
homomorphism $u$ by completion.

\end{document}